\documentclass[a4paper,12pt]{article}
\usepackage{cmap}                        % Поддержка поиска русских слов в PDF (pdflatex)
\usepackage[cp1251]{inputenc}            % Выбор языка и кодировки
\usepackage[english]{babel}
\usepackage[left=3cm,right=2cm,top=2cm,bottom=2cm]{geometry} % поля страницы
\usepackage[ruled,vlined]{algorithm2e}
\usepackage{amssymb}
\usepackage{amsmath, amsthm}
\theoremstyle{plain}
\newtheorem{thm}{Theorem}
\newtheorem{lem}{Lemma}

\newtheorem{cor}{Corollary}[thm]

\theoremstyle{definition}

\begin{document}

\begin{center}\large
\textbf{Constructing regular saturated formations of finite soluble groups\footnote{This work was financially supported by the Belarussian Republican Foundation for Fundamental
Research (BRFFR-RSF M, project F23RNFM-63).}}\normalsize

\smallskip
Viachaslau I. Murashka

 \{mvimath@yandex.ru\}

Department of  Mathematics and Technologies of Programming,

 Francisk Skorina Gomel State University, Gomel, Belarus\end{center}

\begin{abstract}
For a formation $\mathfrak{F}$ of finite groups consider a graph whose vertices are elements of a finite group and two vertices are connected by an edge if and only if they generates non-$\mathfrak{F}$-group as elements of a group. A hereditary formation $\mathfrak{F}$ is called regular if the set of all isolated vertices of the described graph coincides with the intersection of all maximal $\mathfrak{F}$-subgroups in every group. The constructive description of saturated regular formations of soluble groups which improves the results of Lucchini and Nemmi is obtained in this paper. In particular, it is showed that saturated regular non-empty formations of soluble groups are just hereditary formations $\mathfrak{F}$ of soluble groups that contains every group all whose  cyclic primary subgroups are $K$-$\mathfrak{F}$-subnormal. Also we prove that the lattice of saturated regular formations of soluble groups is lattice isomorphic to the Steinitz's lattice.

\textbf{Keywords}: finite group; non-$\mathfrak{F}$-graph of a group,  maximal $\mathfrak{F}$-subgroup, regular formation, hereditary saturated formation, $K$-$\mathfrak{F}$-subnormal subgroup,  Steinitz's lattice.

\textbf{MSC}:  20D10; 06B23.
\end{abstract}

\section*{Introduction}

All groups considered here are finite. One of the classic methods of studying the structure of a group is the graph method, which goes back to the works of Cayley in 1878. The essence of this method is to establish a connection between the characteristics of a group and the graph associated with it. One of the directions of this method is the study of group graphs whose vertices  are the elements of a group. Examples of such graphs are: the Cayley graph, the non-commuting graph \cite{Neumann1976APO}, which was brought to attention by Paul Erd\"os, the non-nilpotent graph \cite{Abdollahi2010}, and many others.

Among such graph's constructions there are important constructions of $\mathfrak{F}$-graphs (see \cite{Delizia2007}) and non-$\mathfrak{F}$-graphs (see \cite{Lucchini2021}). Their main idea is that two elements $x$ and $y$ are connected by an edge if $\langle x, y\rangle\in\mathfrak{F}$ for $\mathfrak{F}$-graphs and $\langle x, y\rangle\not\in\mathfrak{F}$ for non-$\mathfrak{F}$-graphs. Hence the cyclic \cite{Ma2013}, the non-cyclic \cite{Abdollahi2007}, the commuting  \cite{Segev1999}, the non-commuting \cite{Neumann1976APO}, the nilpotent \cite{Das2015} and the non-nilpotent \cite{Abdollahi2010} graphs are particular cases of this construction.

For the study of such graphs the following subset of a group $G$ is important
$$\mathcal{I}_\mathfrak{F}(G)=\{x\mid \langle x, y\rangle\in\mathfrak{F}\,\,\forall y\in G\}.$$
It is easy to see that if $\mathfrak{F}=\mathfrak{A}$ is the class of all abelian groups, then $\mathcal{I}_\mathfrak{A}(G)=\mathrm{Z}(G)$ is the center of $G$. From \cite{Abdollahi2010} it follows that $\mathcal{I}_\mathfrak{N}(G)=\mathrm{Z}_\infty(G)$ is the hypercenter of $G$ where $\mathfrak{N}$ denotes the class of all nilpotent groups. According to \cite{Guralnick2006} $\mathcal{I}_\mathfrak{S}(G)=R(G)$ is the (soluble) radical  of $G$ for the class $\mathfrak{S}$ of all soluble groups.

 It is well known that the intersection of all maximal abelian subgroups is the center of a group. Baer \cite{Baer1953} proved that the intersection of maximal nilpotent subgroups is the hypercenter of a group. It is straightforward to check that the intersection of maximal soluble subgroups is the radical of a group. Recall \cite{h4} that the intersection of all $\mathfrak{F}$-maximal subgroup of a group $G$ is denoted by $\mathrm{Int}_\mathfrak{F}(G)$. According to \cite{Lucchini2021} a formation $\mathfrak{F}$ is called \textbf{\emph{regular}}, if it is hereditary and $\mathrm{Int}_\mathfrak{F}(G)=\mathcal{I}_\mathfrak{F}(G)$ for every group $G$.

 In  \cite{Lucchini2021} all saturated regular formations of soluble groups were described in terms of their  minimal non-$\mathfrak{F}$-groups (or critical groups). The aim of this paper is to improve that description and present the new ways of constructing
  saturated regular formations of soluble groups.

\section{Main results}

Let $\mathfrak{H}$ be a homomorph.  Recall (see \cite[1.2.13]{Kamornikov2003} or \cite[Definition 6.1.4]{s9}) that a subgroup $H$ of  $G$ is called $K$-$\mathfrak{H}$-\emph{subnormal} in $G$ if there is a chain of subgroups
$ H=H_0\subseteq H_1\subseteq\dots\subseteq H_n=G$
with $H_{i-1}\trianglelefteq H_i$ or $H_{i}/\mathrm{Core}_{H_{i}}(H_{i-1})\in\mathfrak{H}$ for all $i=1,\dots,n$.
The class $v^*\mathfrak{H}$ of all groups all whose cyclic primary subgroups are $K$-$\mathfrak{H}$-subnormal was studied in \cite{Murashka2014}. Here we prove

\begin{thm}\label{thmc2}
   If $\mathfrak{H}$ is a hereditary homomorph of soluble groups, then $v^*\mathfrak{H}$ is the least by inclusion saturated regular formation of soluble groups which contains $\mathfrak{H}$. In particular, a hereditary non-empty saturated formation $\mathfrak{F}$ of soluble groups is regular iff $v^*\mathfrak{F}=\mathfrak{F}$.
\end{thm}

Recall \cite[Definition 1]{Vasilev2010} that a subgroup $H$ of $G$ is called $\mathbb{P}$-\emph{subnormal} if $H=G$ or there is a chain of subgroups
$H=H_0< H_1<\dots< H_n=G$ such that $|H_i:H_{i-1}|$ is a prime for $1\leq i\leq n$. The class $v\mathfrak{U}$ of all groups whose cyclic primary subgroups are $\mathbb{P}$-subnormal was studied in \cite{Monakhov2013}. In particular, $v\mathfrak{U}=v^*\mathfrak{U}$ by \cite[Theorem B]{Murashka2014}.

\begin{cor}
  The class $v\mathfrak{U}$ is a regular formation.
\end{cor}

Recall \cite{Semenchuck1984} that a formation  $\mathfrak{F}$ is said to have  the Shemetkov property in $\mathfrak{S}$ if every soluble minimal non-$\mathfrak{F}$-group is either a Schmidt group or a group of prime order. According to \cite{Skiba1990} if a hereditary  formation  $\mathfrak{F}$ of soluble groups  have  the Shemetkov property in $\mathfrak{S}$, then it is saturated. Hence   from \cite[Corollary E.2]{Murashka2014} and Theorem \ref{thmc2} it follows

\begin{cor}
 Let $\mathfrak{F}$ be a hereditary formation of soluble groups     with the Shemetkov property in $\mathfrak{S}$ and $\pi(\mathfrak{F})=\mathbb{P}$.  Then  $\mathfrak{F}$ is  regular.
\end{cor}

It is known \cite[Proposition 2]{Ito1951} that a non-$p$-nilpotent group all whose proper subgroups are $p$-nilpotent is a Schmiodt group.

\begin{cor}
 Let $p$ be a prime. The class of all soluble $p$-nilpotent groups is a regular formation.
\end{cor}

Let $\varphi$ be some linear ordering on $\mathbb{P}$. Recall \cite[IV, Examples 3.4(g)]{s8} that a group $G$ is called a Sylow tower group of type $\varphi$  if  it has normal Hall $\{p_1, \dots, p_t\}$-subgroups for all $1\leq t\leq k$ where $\pi(G)=\{p_1,\dots, p_k\}$ and $p_i>_{\varphi}p_j$ for $i<j$. It is clear that $\varphi$-dispersive group is soluble. Using \cite[Proposition 2]{Ito1951} one can show that every minimal non $\varphi$-dispersive group is a Schmidt group.

\begin{cor}
 Let $\varphi$ be some linear ordering on $\mathbb{P}$. Then the class of all Sylow tower groups of type  $\varphi$ is a regular formation.
 %In particular, the class of groups with the Sylow tower property is a regular formation.
\end{cor}

Let $\mathfrak{F}_1$ and $\mathfrak{F}_2$ be  saturated regular non-empty formations of soluble groups. Denote by $\mathfrak{F}_1\vee_{reg}\mathfrak{F}_2$ and $\mathfrak{F}_1\wedge_{reg}\mathfrak{F}_2$ the least by inclusion saturated regular formation which contains $\mathfrak{F}_1$ and $\mathfrak{F}_2$, and  the greatest by inclusion saturated regular formation contained in $\mathfrak{F}_1$ and $\mathfrak{F}_2$   respectively. In the proof of Theorem \ref{thmc2} we will show that these operations are well defined.
\begin{cor}
  Let $\mathfrak{F}_1$ and $\mathfrak{F}_2$ be  saturated regular non-empty formations of soluble groups. Then   $\mathfrak{F}_1\wedge_{reg}\mathfrak{F}_2=\mathfrak{F}_1\cap  \mathfrak{F}_2$  and $\mathfrak{F}_1\vee_{reg}\mathfrak{F}_2=v^*(\mathfrak{F}_1\cup  \mathfrak{F}_2)$.
\end{cor}

The proof of Theorem \ref{thmc2} is based on the fact, that the values of the closure operation $v^*$ on hereditary homomorphs are solubly saturated formations.

\begin{thm}\label{thmc3}
 Let $\mathfrak{H}$ be a hereditary  homomorh. Then $v^*\mathfrak{H}$ is a solubly saturated formation. In particular, if $\mathfrak{H}$ consists of soluble groups, then $v^*\mathfrak{H}$ is a saturated formation.
\end{thm}

Recall that a supernatural number is the following formal product $\omega=\prod_{p\in\mathbb{P}}p^{v_p(\omega)}$ where $v_p(\omega)$ can be zero, a natural number or infinity. We use $\mathbb{SN}$ to denote the set of all supernatural numbers. If there are no infinities among $v_p(\omega)$ and only a finite number of $v_p(\omega)$ are non-zero, then $\omega$ is just a natural number. Recall that $\omega_1$ divides $\omega_2$ if $v_p(\omega_1)\leq v_p(\omega_2)$ for all $p\in\mathbb{P}$.
For a class of groups $\mathfrak{X}$ and a supernatural number $\omega$ we use $\mathfrak{X}(\omega)$ to denote the class of all $\mathfrak{X}$-groups of exponent dividing $\omega$.

\begin{thm}\label{thmc1}
The following statements hold:

$(1)$ Let $f: \mathbb{P}\mapsto\mathbb{SN}$ be a function with $v_p(f(p))=\infty$ for all $p\in\mathbb{P}$. Then
 $$\bigcap_{p\in\mathbb{P}}\mathfrak{S}_{p'}\mathfrak{S}(f(p))$$
 is a  saturated regular  formation of soluble groups.

 $(2)$ If  $\mathfrak{F}$ is a saturated regular  non-empty formation of soluble groups, then there exists a function $f: \mathbb{P}\mapsto\mathbb{SN}$ with $v_p(f(p))=\infty$ for all $p\in\mathbb{P}$ such that
 $$\mathfrak{F}=\bigcap_{p\in\mathbb{P}}\mathfrak{S}_{p'}\mathfrak{S}(f(p)).$$
\end{thm}

For supernatural numbers $\omega_1$ and $\omega_2$ recall that $$lcm(\omega_1,\omega_2)=\prod_{p\in\mathbb{P}}p^{\max\{v_p(\omega_1),v_p(\omega_2)\}}\textrm{ and }gcd(\omega_1,\omega_2)=\prod_{p\in\mathbb{P}}p^{\min\{v_p(\omega_1),v_p(\omega_2)\}}.$$
The set $\mathbb{SN}$ with operations $lcm$ and $gcd$ is a complete distributive lattice which is called the Steinitz's lattice \cite{Steinitz1910}, \cite[Definition 2]{Bezushchak2016}.
Denote by $\mathrm{Reg}$ the set of all  saturated regular non-empty formations of soluble groups.

\begin{thm}\label{thmc4}
 The lattice $(\mathrm{Reg}, \vee_{reg}, \wedge_{reg})$ is lattice isomorphic to the Steinitz's lattice. In particular it is a complete distributive lattice.
\end{thm}

According to \cite[Definition 3]{Bezushchak2016} a supernatural number $\omega$ is called complete if  $v_p(\omega)\in\{0,\infty\}$ for all $p\in\mathbb{P}$. Denote by $\mathrm{Shem}$ the set of all hereditary formations of soluble groups with the Shemetkov property in $\mathfrak{S}$ which contains $\mathfrak{N}$.

\begin{cor}\label{cor5}
 The lattice $(\mathrm{Shem}, \vee_{reg}, \wedge_{reg})$ is lattice isomorphic to the sublattice of all complete supernatural numbers of Steinitz's lattice. In particular it is a complete distributive and complemented lattice.
\end{cor}

From this corollary the main result of \cite{BallesterBolinches2024} follows

\begin{cor}[{\cite[Theorem A(1, 3)]{BallesterBolinches2024}}]\label{cor6}
The set $\mathrm{Shem}$ with natural partial order is a complete, distributive and complemented lattice with the least and the greatest elements.\end{cor}

\section{Preliminaries}

All unexplained notations and terminologies are standard. The reader is referred to \cite{s9, s8, s6} if necessary.
Recall that $\Phi(G)$ is the Frattini subgroup of $G$; $\mathrm{Soc}(G)$ denotes the socle of $G$; a generalization of the Fitting subgroup $\mathrm{\tilde F}(G)$ of a group $G$ is defined by $\Phi(G)\subseteq \mathrm{\tilde F}(G)$ and $\mathrm{\tilde F}(G)/\Phi(G)=\mathrm{Soc}(G/\Phi(G))$; $\mathbb{P}$ denotes the set of all primes; $\pi(G)$ is the set of all prime divisors of  $|G|$;  $\pi(\mathfrak{X})=\underset{G\in\mathfrak{X}}\cup\pi(G)$.

According to \cite[B, Theorem 12.4]{s8} there exists the uniquely-determined group $G$, denoted by $E(n|p)$, such that
$G$ has an abelian normal subgroup $A$  of exponent $p$,    $A$ has a complement $C\simeq   Z_n$ in $G$, and    $C$ acts faithfully and indecomposably on $A$.

Recall that a \emph{homomorph} $\mathfrak{H}$ is a class of groups which is closed under taking epimorphic images (i.e. from $G\in\mathfrak{H}$ and $N\trianglelefteq G$ it follows that $G/N\in\mathfrak{H}$); a homomorph $\mathfrak{F}$ is called \emph{formation} if it is closed under taking  subdirect products (i.e. from $G/N_1\in\mathfrak{F}$ and $G/N_2\in\mathfrak{F}$ it follows that $G/(N_1\cap N_2)\in\mathfrak{F}$). A formation $\mathfrak{F}$ is said to be: \emph{saturated}  if $G\in\mathfrak{F}$
whenever $G/\Phi(G)\in\mathfrak{F}$; \emph{solubly saturated}  if $G\in\mathfrak{F}$
whenever $G/\Phi(R(G))\in\mathfrak{F}$; \emph{hereditary} if $H\in \mathfrak{F}$ whenever $H\leq G\in \mathfrak{F}$.

  Let $\mathfrak{X}$ be a class of groups.  A chief factor $H/K$ of  $G$ is called   $\mathfrak{X}$-\emph{central} in $G$ provided    $(H/K)\rtimes (G/C_G(H/K))\in\mathfrak{X}$ (see \cite[p. 127--128]{s6}). A normal subgroup $N$ of $G$ is said to be $\mathfrak{X}$-\emph{hypercentral} in $G$ if $N=1$ or $N\neq 1$ and every chief factor of $G$ below $N$ is $\mathfrak{X}$-central. The symbol $\mathrm{Z}_\mathfrak{X}(G)$ denotes the $\mathfrak{X}$-\emph{hypercenter} of $G$, that is, the greatest normal $\mathfrak{X}$-hypercentral subgroup of $G$ (see  \cite[Lemma 14.1]{s6}).

  A formation $\mathfrak{F}$ is called $Z$-\emph{saturated} \cite{Murashka2022a} if $\mathfrak{F}=(G\mid G=\mathrm{Z}_\mathfrak{F}(G))$.

\begin{lem}\label{lem1}
  Let $\mathfrak{H}$ be a hereditary homomorph. Then

  $(1)$ $v^*\mathfrak{H}$ is a hereditary formation \cite[Theorem A(3)]{Murashka2014}.

  $(2)$ $v^*\mathfrak{H}=v^*(v^*\mathfrak{H})$ \cite[Theorem A(4)]{Murashka2014}.

  $(3)$ $v^*\mathfrak{H}$ is a $Z$-saturated   formation \cite[Proposition 6]{Murashka2022}.
\end{lem}

Recall that $\mathbb{F}_p$ denotes a field with $p$ elements for a prime $p$.
Whenever
$V$ is a $G$-module over $\mathbb{F}_p$, $V\rtimes G$ denotes the semidirect product of $V$ with $G$
corresponding to the action of $G$ on $V$ as $G$-module.
A formation $\mathfrak{F}$ is said  to satisfy \emph{the tensor product property in the class of all
groups} \cite[Definition 1(1)]{Murashka2022a} if for any $p\in\mathbb{P}$ the following condition holds for any two   $G$-modules $V$ and $W$ over $\mathbb{F}_p$:
$$V \rtimes G \in \mathfrak{F} \textrm{ and } W \rtimes G \in\mathfrak{F} \Rightarrow (V\otimes W) \rtimes G \in \mathfrak{F}.$$
The similar condition in the class of all soluble groups was considered in \cite{BallesterBolinches1999}.

Let $\mathfrak{X}$ be a class of groups. A group $G\not\in\mathfrak{X}$ is called \emph{a minimal non-$\mathfrak{X}$-group} if all its proper subgroups belong $\mathfrak{X}$. A minimal non-$\mathfrak{X}$-group  is called \emph{strongly critical for $\mathfrak{X}$} if every its proper quotient belongs $\mathfrak{X}$.
\begin{thm}[Lucchini and Nemmi \cite{Lucchini2021}]\label{lnt}
  Let $\mathfrak{F}$ be a hereditary saturated formation, with $\mathfrak{A} \subseteq \mathfrak{F} \subseteq \mathfrak{S}$.
Then $\mathfrak{F}$ is regular if and only if every finite group $G$ which is soluble and
strongly critical for $\mathfrak{F}$ has the property that $G/\mathrm{Soc}(G)$ is cyclic.
\end{thm}

We use $exp(G)$ to denote the \emph{exponent} of $G$, i.e. the least common multiple of the orders of all elements of the group. For an infinite set of natural numbers $A$ let $\max A=\infty$. For a class of groups $\mathfrak{X}$ let $$exp(\mathfrak{X})=lcm\{exp(G)\mid G\in \mathfrak{X}\}=\prod_{p\in\mathbb{P}}p^{\max\{v_p(G)\mid G\in\mathfrak{X}\}}.$$
Note that  $exp(\mathfrak{X})$ is a supernatural number.

Let $\mathfrak{F}$ be a non-empty formation. Then is every group $G$ exists the least normal subgroup $G^\mathfrak{F}$ with $G/G^\mathfrak{F}\in\mathfrak{F}$. If $\mathfrak{H}$ is a formation, then $\mathfrak{HF}=(G\mid G^\mathfrak{F}\in\mathfrak{H})$ is also formation \cite[IV, Theorem 1.8(a)]{s8}.

A function of the form $f: \mathbb{P}\rightarrow\{formations\}$ is called a \emph{formation function}. Recall \cite[IV, Definitions 3.1]{s8} that a formation $\mathfrak{F}$ is called \emph{local} if $\mathfrak{F}=(G\,|\, G/C_G(H/K)\in f(p)$ for every $p\in\pi(H/K)$ and every chief factor $H/K$ of $G$) for some formation function $f$. In this case $f$ is called a \emph{local definition} of $\mathfrak{F}$.
Recall \cite[IV, Proposition 3.8]{s8} that if $\mathfrak{F}$ is a local formation, there exists an unique formation function $F$, defining $\mathfrak{F}$, such that    $F(p) = \mathfrak{N}_p F(p)\subseteq\mathfrak{F}$ for every $p\in\mathbb{P}$. In this case $F$ is called the \emph{canonical local definition}  of $\mathfrak{F}$.

\section{Proves of the main results}

 \subsection{Proof of Theorem \ref{thmc3}}

$(a)$  \emph{$\mathfrak{F}=v^*\mathfrak{H}=v^*\mathfrak{F}$ is a hereditary $Z$-saturated formation.}

Let $\mathfrak{F}=v^*\mathfrak{H}$.
From  Lemma \ref{lem1} it follows that $\mathfrak{F}=v^*\mathfrak{H}=v^*(v^*\mathfrak{H})=v^*\mathfrak{F}$ is a hereditary $Z$-saturated formation.

 $(b)$ \emph{Let $p\neq q$ be primes. If $E(q^k|p)\in\mathfrak{F}$, then $E(q^t|p)\in\mathfrak{F}$ for all $t\leq k$.}

Let $E(q^k|p)\simeq G=NC$ where $N$ is an elementary abelian $p$-group and $C$ is a cyclic group of order $q^k$. Note that $C_G(N)=N$. Let $C_1$ be a cyclic subgroup of order $q^t$ for $0\leq t<k$. Note that if $t=0$, then $E(q^0|p)\simeq Z_p\in\mathfrak{F}$. Hence we may assume that $t>0$.  From $gcd(|N|,|C_1|)=1$ and \cite[A, Theorem 11.6]{s8}    it follows that $N$ is the direct product of minimal normal subgroups $N_i$ of $NC_1$, $1\leq i\leq n$. Let $C_2\simeq Z_q$ be a subgroup of $C_1$. Note that $C_2$ is the unique minimal normal subgroup of $C_1$. Assume that $C_{C_1}(N_i)\neq 1$ for all $1\leq i\leq n$. Then $$1\neq C_2\leq\bigcap_{i=1}^n  C_{C_1}(N_i)=C_{C_1}(N)\leq C_C(N)=1,$$ a contradiction. Thus there is $j$ with $C_{C_1}(N_j)= 1$. Since $\mathfrak{F}$ is hereditary, $N_jC_1\in\mathfrak{F}$. Now $N_jC_1$ has the unique minimal normal subgroup $N_j$, $N_j$ is elementary   abelian $p$-subgroup, $C_1\simeq Z_{q^t}$ is the complement of $N_j$ and $C_1$ acts faithfully and irreducibly on $N_j$. Hence $N_jC_1\simeq E(q^t|p)$.  Thus $E(q^t|p)\in\mathfrak{F}$ for all $t\leq k$.

 $(c)$ \emph{Let $N$ and $M$ be $C$-modules over $\mathbb{F}_p$ with $C_C(N)\leq C_C(M)$ where $p\neq q$ are primes and $C\simeq Z_{q^l}$ for some $l\geq 1$. If  $N\rtimes C\in\mathfrak{F}$, then $M\rtimes C\in\mathfrak{F}$}.

Note that $C_C(N)\trianglelefteq NC$.
By analogy with the proof of $(b)$, $N$ is the direct product of minimal normal subgroups $N_i$ of $NC$, $1\leq i\leq n$, and there is $j$ with $C_C(N_j)=C_C(N)$. Let $C/C_C(N)\simeq Z_{q^k}$. If $k=0$, then $M\rtimes C\in\mathfrak{F}$ as the direct product of $\mathfrak{F}$-groups of the form $Z_{q^l}$ and $Z_p$. Assume now that $k>0$. Note that $N_jC/C_C(N)\simeq E(q^k|p)$. Now $E(q^t|p)\in \mathfrak{F}$ for all $t\leq k$ by $(b)$.

Let $H/K$ be a chief factor of $G= M\rtimes C$ below $M$. Then $MC_C(N)\leq MC_C(M)\leq C_G(H/K)$. Hence $(H/K)\rtimes G/C_G(H/K)\simeq E(q^t|p)\in\mathfrak{F}$ for some $t\leq k$. If $H/K$ is a chief factor of $G$ above $N$, then $(H/K)\rtimes G/C_G(H/K)\simeq Z_q\in\mathfrak{F}$. Thus $G=\mathrm{Z}_\mathfrak{F}(G)\in\mathfrak{F}$.

$(d)$ \emph{$\mathfrak{F}$ has the tensor product property in the class of all groups}.

Let $V$ and $M$ be  $G$-modules over $\mathbb{F}_p$ with $V\rtimes G, M\rtimes G\in\mathfrak{F}$, $T=(V\otimes M)\rtimes G$ and $N=V\otimes M$. Since $\mathfrak{F}$ is a hereditary formation, $G\in\mathfrak{F}$ and $Z_p\in\mathfrak{F}$.

$(d.1)$ \emph{$NC\in \mathfrak{F}$ for every cyclic primary subgroup $C$ of $G$}.

From $\mathfrak{F}=v^*\mathfrak{F}$  it follows that $\mathfrak{F}$ contains all $p$-groups. Therefore if $C$ is a $p$-group, then  $NC\in\mathfrak{F}$. Assume that $C$ is  a $q$-group for some prime $q\neq p$.

According to \cite[B, Lemma 1.1]{s6} the action of $G$ on pure tensors $v\otimes m$ is defined by $(v\otimes m)^g=v^g\otimes m^g$. Hence if $g\in C_G(V)\cap C_G(M)$, then $g\in C_G(V\otimes M)$. Let $C_1=C_C(V), C_2=C_C(M)$ and $C_3=C_C(N)$. Since $C$ has the unique chief series, $C_1\cap C_2\in\{C_1, C_2\}$. WLOG assume that $C_1\cap C_2=C_1$. Now $C_1\leq C_3$. Note that $V$ and $N$ are  $C$-modules over $\mathbb{F}_p$ with $C_C(V)\leq C_C(N)$ and $V\rtimes C\in \mathfrak{F}$ (as a subgroup of an $\mathfrak{F}$-group $V\rtimes G$). Thus $NC\in\mathfrak{F}$ by $(c)$.

$(d.2)$ \emph{If $Z$ is a cyclic primary subgroup of $T$, then there is a cyclic primary subgroup $C$ of $G$ with $NC=NZ$}.

We have that $NZ=NZ\cap T=NZ\cap NG=N(NZ\cap G)$. From $$NZ/N\simeq  N(NZ\cap G)/N\simeq  (NZ\cap G)/(NZ\cap G\cap N)\simeq(NZ\cap G)$$ it follows that $C=NZ\cap G$ is a cyclic primary subgroup of $G$ with $NC=NZ$.

$(d.3)$ \emph{All cyclic primary subgroups of $T$ are $K$-$\mathfrak{F}$-subnormal, i.e. $T\in\mathfrak{F}$}.

Let $Z$ be a cyclic primary subgroup of $T$. Then there is a cyclic primary subgroup $C$ of $G$ with $NC=NZ$ by $(d.2)$. Hence $NC\in \mathfrak{F}$ by $(d.1)$. It means that $Z$ is $K$-$\mathfrak{F}$-subnormal in $NC$. From $G\in \mathfrak{F}$ it follows that $CN/N\simeq C$ is $K$-$\mathfrak{F}$-subnormal in $T/N\simeq G$. Hence $CN$ is $K$-$\mathfrak{F}$-subnormal in $T$ by \cite[Lemma 6.1.6(2)]{s9}. Thus $Z$ is $K$-$\mathfrak{F}$-subnormal in $T$ by \cite[Lemma 6.1.6(1)]{s9}. Therefore $T\in v^*\mathfrak{F}=\mathfrak{F}$.

It means that $\mathfrak{F}$ has the tensor product property in the class of all groups.

$(e)$ \emph{$\mathfrak{F}$ is solubly saturated. Moreover if $\mathfrak{H}$ consists of soluble groups, then $\mathfrak{F}$ is a saturated formation}.

Since a $Z$-saturated formation $\mathfrak{F}$ has the tensor product property in the class of all groups, it is solubly saturated by \cite[Theorem 4]{Murashka2022a}. If $\mathfrak{H}$ consists of soluble groups, then $\mathfrak{F}$ also consists of soluble groups. It is well known that solubly saturated formations of soluble groups are saturated.

 \subsection{Proof of Theorem \ref{thmc2}}

Suppose that $\mathfrak{F}$ is a saturated regular non-empty formation of soluble groups.

From $1\in\mathrm{Int}_\mathfrak{F}(G)=\mathcal{I}_\mathfrak{F}(G)$ for every group $G$ it follows that every cyclic group belongs $\mathfrak{F}$. Hence $\pi(\mathfrak{F})=\mathbb{P}$.

Assume that $G$ is a nonsoluble  minimal non-$\mathfrak{F}$-group. Then $G$ is a minimal nonsoluble group. Note that $G/\Phi(G)$ is a minimal simple non-abelian group. It means that $G=\mathrm{\tilde F}(G)$ and $G/\mathrm{\tilde F}(G)\simeq 1$ is a cyclic primary group.

Let $G$ be a soluble minimal non-$\mathfrak{F}$-group. Since $\mathfrak{F}$ is saturated,  $\overline{G}=G/\Phi(G)$ is also a soluble minimal non-$\mathfrak{F}$-group.  By Gasch\"utz-Lubeseder-Schmid Theorem \cite[IV, Theorem 4.6]{s8}  $\mathfrak{F}$ is local. Let $F$ be the canonical local definition of $\mathfrak{F}$. Then $\overline{G}\simeq N\rtimes M$ where $N=\mathrm{Soc}(\overline{G})$ is a $p$-group and $M$ is a minimal non-$F(p)$-group by \cite[Theorem 1.2]{Semenchuk1979}. Note that every proper quotient of $\overline{G}$ is an $\mathfrak{F}$-group, i.e. $G$ is strongly critical for $\mathfrak{F}$.  Since $\mathfrak{F}$ is regular, $M$ is cyclic by Theorem \ref{lnt}. Assume that $M$ is not primary. Then every its Sylow subgroup  belongs $F(p)$. Hence $M\in F(p)$ as the direct product of its Sylow subgroups, a contradiction. Thus $G/\mathrm{\tilde F}(G)\simeq \overline{G}/\mathrm{Soc}(\overline{G})$ is a cyclic primary group for every soluble   minimal non-$\mathfrak{F}$-group $G$.

 Therefore $G/\mathrm{\tilde F}(G)$ is a cyclic primary group for every    minimal non-$\mathfrak{F}$-group $G$ and $\pi(\mathfrak{F})=\mathbb{P}$. Thus $\mathfrak{F}=v^*\mathfrak{F}$ by \cite[Corollary E.2]{Murashka2014}.

  Suppose that $\mathfrak{H}$ is a hereditary homomorph of soluble groups. From Theorem  \ref{thmc3} it follows that $\mathfrak{F}=v^*\mathfrak{H}$ is a hereditary saturated formation. According to (2) of Lemma  \ref{lem1} we have that $v^*\mathfrak{F}=v^*(v^*\mathfrak{H})=v^*\mathfrak{H}=\mathfrak{F}$. Now according to \cite[Corollary E.2]{Murashka2014} $G/\mathrm{\tilde F}(G)$ is a cyclic primary group for every minimal non-$\mathfrak{H}$-group $G$. Since $\mathfrak{F}$ is saturated, a soluble  minimal non-$\mathfrak{F}$-group is strongly critical for $\mathfrak{F}$ iff $\Phi(G)=1$. In this case  $G/\mathrm{\tilde F}(G)=G/\mathrm{Soc}(G)$   is a cyclic (primary) group. Thus $\mathfrak{H}$ is regular by Theorem \ref{lnt}. So if $\mathfrak{H}$ is a hereditary homomorph, then $v^*\mathfrak{H}$ is a saturated regular formation.

Let $\{\mathfrak{F}_i\mid i\in I\}$  be a set of saturated non-empty regular formations and $\mathfrak{F}=\cap_{i\in I}\mathfrak{F}_i$.  Then $\mathfrak{F}$ is a hereditary saturated formation. Note that $\mathfrak{N}\subseteq \mathfrak{F}$. From the definition of   $\mathfrak{F}$ it follows that every strongly critical for $\mathfrak{F}$ group $G$ is  strongly critical for $\mathfrak{F}_i$ group for some $i\in I$.  Then $G/\mathrm{Soc}(G)$ is cyclic by Theorem \ref{lnt}. Now $\mathfrak{F}$ is regular by Theorem \ref{lnt}. Hence  $\mathfrak{F}$ is the greatest by inclusion saturated regular formation contained in $\mathfrak{F}_i$ for all $i\in I$. Note that for a class $\mathfrak{X}$ the intersection of all non-empty saturated regular formations which contains $\mathfrak{X}$ is the least by inclusion saturated regular formation which contains $\mathfrak{X}$.

  Let $\mathfrak{H}$ be a hereditary homomorph and $\mathfrak{F}$ be the least by inclusion saturated regular formation with $\mathfrak{H}\subseteq\mathfrak{F}$. Note that $v^*\mathfrak{F}=\mathfrak{F}$. Now $v^*\mathfrak{H}\subseteq v^*\mathfrak{F}=\mathfrak{F}$  by \cite[Theorem A(5)]{Murashka2014} and $ v^*\mathfrak{H}$ is a saturated regular formation. Thus $v^*\mathfrak{H}=\mathfrak{F}$.

 \subsection{Proof of Theorem \ref{thmc1}}

Let $\omega\in \mathbb{SN}$. Note that a soluble group $G\in\mathfrak{S}(\omega)$ if and only if every cyclic primary subgroup of $G$ is in $\mathfrak{S}(\omega)$. Also note that if $v_p(f(p))=\infty$, then $\mathfrak{S}(f(p))=\mathfrak{N}_p\mathfrak{S}(f(p))$ because  $p'$-groups belonging to $\mathfrak{S}(f(p))$ and $\mathfrak{N}_p\mathfrak{S}(f(p))$ are the same and all $p$-groups belong to these classes.

(1) Let $\mathfrak{H}=\bigcap_{p\in\mathbb{P}}\mathfrak{G}_{p'}\mathfrak{S}(f(p))=\bigcap_{p\in\mathbb{P}}\mathfrak{G}_{p'}\mathfrak{N}_p\mathfrak{S}(f(p))$. Then $\mathfrak{F}$ is a local formation locally defined by $h(p)=\mathfrak{S}(f(p))$ by \cite[IV, Theorem 3.2]{s8}. Hence $\mathfrak{H}$ consists of soluble groups. Now  $\mathfrak{H}=\bigcap_{p\in\mathbb{P}}\mathfrak{S}_{p'}\mathfrak{S}(f(p))$. By Gasch\"utz-Lubeseder-Schmid Theorem \cite[IV, Theorem 4.6]{s8}  $\mathfrak{F}$ is saturated.

Note that $h(p)=\mathfrak{S}(f(p))$ is hereditary for all $p\in\mathbb{P}$. Now $\mathfrak{H}$ is hereditary too by \cite[IV, Proposition 3.14]{s8}. Assume that $v^*\mathfrak{H}\neq\mathfrak{H}$. Hence there is a group $G\in v^*\mathfrak{H}\setminus\mathfrak{H}$. We may assume that $G$ is the minimal group with this property. Since $v^*\mathfrak{H}$ and $\mathfrak{H}$ are formations of soluble groups and $\mathfrak{H}$ is saturated, we see that $G$ is a minimal non-$\mathfrak{H}$-group and soluble, $\Phi(G)=1$ and $G$ has the unique minimal normal subgroup $N$. Hence $G=N\rtimes M$ and $M\in \mathfrak{H}$. Note that all $p$-subgroups of $M$ belong $h(p)$.

Let $C$ be a cyclic primary $p'$-subgroup of $M$.
From $G\in v^*\mathfrak{H}\setminus\mathfrak{H}$ it follows that $NC$ is a proper subgroup of $G$. Hence $NC\in\mathfrak{H}$. From $C_G(N)=N$ it follows that $\mathrm{O}_{p',p}(NC)=N$. Hence $C\simeq NC/N=NC/\mathrm{O}_{p',p}(NC)\in h(p)$.  It means that $M\simeq G/C_G(N)\in (h(p)\cap\mathfrak{H})$. Therefore $G\in\mathfrak{H}$, the contradiction. Thus $v^*\mathfrak{H}=\mathfrak{H}$. Now $\mathfrak{H}$ is regular by Theorem \ref{thmc2}.

(2) Let $\mathfrak{F}$ be a saturated regular non-empty formation of soluble groups. Then $v^*\mathfrak{F}=\mathfrak{F}$ by
Theorem \ref{thmc2}. Since $\mathfrak{F}\neq\emptyset$, we see that $\mathfrak{F}$ is local by Gasch\"utz-Lubeseder-Schmid Theorem.

Let $F$ be the canonical local definition of $\mathfrak{F}$. Then the canonical local definition of $v^*\mathfrak{F}$ is the class of all $v^*\mathfrak{F}$-groups all whose cyclic primary subgroups belong to $F(p)$ by \cite[Theorems B and D]{Murashka2014}. It means that in our case  $F(p)=\mathfrak{F}(exp(F(p)))$. Let $f(p)=exp(F(p))$. From $\mathfrak{N}_p\subseteq F(p)$ for all $p\in\mathbb{P}$ it follows that $v_p(f(p))=\infty$ for all $p\in\mathbb{P}$.

Let $G(p)=\mathfrak{S}(f(p))$ and $\mathfrak{H}=LF(G)$. Note that $\mathfrak{H}\subseteq\mathfrak{S}$.
From $F(p)\subseteq G(p)$ for all $p\in\mathbb{P}$ it follows that $\mathfrak{F}\subseteq \mathfrak{H}$. Assume that $\mathfrak{F}\neq \mathfrak{H}$. Let $G$ be a minimal order group from $\mathfrak{H}\setminus\mathfrak{F}$. Since $\mathfrak{F}$ is a hereditary saturated formation, we see that $\Phi(G)=1$, the soluble group $G$ has the unique minimal normal subgroup $N$, $N$ is a $p$-group for some $p\in\mathbb{P}$,  $C_G(N)=N$ and $G/N\in\mathfrak{F}$.

Assume that $NC<G$ for every cyclic primary $p'$-subgroup $C$ of $G$. Then $NC\in\mathfrak{F}$. From $N=C_G(N)$ it follows that $C\simeq NC/N=NC/\mathrm{O}_{p',p}(NC)\in F(p)$.  Note that $\mathfrak{N}_p\subseteq F(p) $. Now all cyclic primary subgroups of an $\mathfrak{F}$-group $G/N$ belong $F(p)$. Thus $G/N=G/C_G(N)\in F(p)$. It means that $G\in\mathfrak{F}$, a contradiction.

Therefore there is a cyclic primary $p'$-subgroup $C$ of $G$ with $G=NC$. Note that $C\in\mathfrak{F}$. Hence $exp(C)$ does not divide $exp(F(p))$.  Now $C\simeq G/C_G(N)\in G(p)=\mathfrak{S}(exp(F(p)))$. Hence the exponent of $C$ divides $exp(F(p))$, the final contradiction. Thus $\mathfrak{F}=\mathfrak{H}$. By analogy with the proof of $(1)$ we see that $\mathfrak{F}=\bigcap_{p\in\mathbb{P}}\mathfrak{S}_{p'}\mathfrak{S}(exp(F(p)))$.

\subsection{Proof of Theorem \ref{thmc4}}

In this subsection $\mathfrak{F}, \mathfrak{F}_1$ and $\mathfrak{F}_2$ denote some elements of $\mathrm{Reg}$. Let $F$, $F_1$ and $F_2$ be the canonical local definitions of $\mathfrak{F}, \mathfrak{F}_1$ and $\mathfrak{F}_2$ respectively.

Denote by $p_i$ the $i$-th prime number (i.e. $p_1=2, p_2=3,\dots$).
Note that the set $\mathcal{N}=\{(n_1, n_2)\mid n_1\neq n_2, n_1, n_2\in\mathbb{N}\}$ is countable. So there is a bijection $f: \mathbb{N}\rightarrow \mathcal{N}$. For $i\in\mathbb{N}$ if $f(i)=(n_1, n_2)$, then denote by  $\pi_1(i)=n_1$ and $\pi_2(i)=n_2$.
Let a function $\mathfrak{f}:\textrm{Reg}\rightarrow\mathbb{SN}$ be defined by
$$\mathfrak{f}(\mathfrak{F})=\prod_{i=1}^\infty p_i^{v_{\pi_2(i)}(exp(F(p_{\pi_1(i)})))}.$$
We claim that $\mathfrak{f}$ is the lattice isomorphism between $(\mathrm{Reg},\vee_{reg}, \wedge_{reg})$ and $(\mathbb{SN}, lcm, gcd)$.

 $(a)$ \emph{$\mathfrak{f}$ is injective}.

%According to $(2)$ of the proof of Theorem \ref{thmc1} .

Assume that $\mathfrak{f}(\mathfrak{F}_1)=\mathfrak{f}(\mathfrak{F}_2)$. Since  $f$ is the bijection between $\mathbb{N}$ and $\mathcal{N}$, we see that $v_q(exp(F_1(p)))=v_q(exp(F_2(p)))$ for all primes $p\neq q$.   From $\mathfrak{N}_p\subseteq F_1(p)\cap F_2(p)$ for all $p\in\mathbb{P}$ it follows that $v_p(exp(F_1(p)))=v_p(exp(F_2(p)))=\infty$ for all $p\in\mathbb{P}$. It means that $exp(F_1(p))=exp(F_2(p))$ for all $p\in\mathbb{P}$. In the proof of Theorem \ref{thmc1} we proved that  $\mathfrak{F}_1$ and $\mathfrak{F}_2$ can be defined by $G_1$ and $G_2$ where $G_i(p)=\mathfrak{S}(exp(F_i(p)))$ for all $i\in\{1,2\}, p\in\mathbb{P}$. Thus $G_1(p)=G_2(p)$ for all $p\in\mathbb{P}$. It means that $\mathfrak{F}_1=\mathfrak{F}_2$. Hence $\mathfrak{f}$ is injective.

 $(b)$ \emph{$\mathfrak{f}$ is surjective}.

 Let $\omega\in\mathbb{SN}$. Denote by $$\alpha(p_k)=p_k^\infty\prod_{j=1, j\neq k}^\infty p_j^{v_{p_{f^{-1}(k,j)}}(\omega)}.$$
Note that $\alpha: \mathbb{P}\rightarrow\mathbb{SN}$ is a function with $v_p(\alpha(p))=\infty$. Let $G(p)=\mathfrak{S}(\alpha(p))$ for all $p\in\mathbb{P}$. Then formation $\mathfrak{H}$, locally defined by $G$,   is a saturated regular formation by Theorem \ref{thmc1}. It is straightforward to check that $\mathfrak{f}(\mathfrak{H})=\omega$.

$(c)$ \emph{If $\mathfrak{F}=\mathfrak{F}_1\vee_{reg}\mathfrak{F}_2$, then $F(p)=\mathfrak{F}(lcm(exp(F_1(p)), exp(F_2(p))))$ for all $p\in\mathbb{P}$. In particular, $\mathfrak{f}(\mathfrak{F}_1\vee_{reg}\mathfrak{F}_2)=lcm(\mathfrak{f}(\mathfrak{F}_1),\mathfrak{f}(\mathfrak{F}_2))$}.

Let $H(p)=\mathfrak{S}(lcm(exp(F_1(p)), exp(F_2(p))))$ for all $p\in\mathbb{P}$. Then formation $\mathfrak{H}$, locally defined by $H$,   is a saturated regular formation by Theorem \ref{thmc1}. Note that $F_1(p)\cup F_2(p)\subseteq F(p)$. It means that $lcm(exp(F_1(p)), exp(F_2(p)))$ divides $exp(F(p))$. According to the proof of Theorem \ref{thmc1} $\mathfrak{F}$ can be locally defined by $G$ where $G(p)=\mathfrak{S}(exp(F(p)))$ for all $p\in\mathbb{P}$. Hence $H(p)\subseteq G(p)$ for all $p\in\mathbb{P}$. It means that $\mathfrak{H}\subseteq \mathfrak{F}$. Since $\mathfrak{F}$ is the least by inclusion element from $\mathrm{Reg}$ which contains $\mathfrak{F}_1$ and $\mathfrak{F}_2$ and $\mathfrak{H}\in\mathrm{Reg}$ also contains $\mathfrak{F}_1$ and $\mathfrak{F}_2$, $\mathfrak{F}=\mathfrak{H}$. Now
\begin{multline*}
  \mathfrak{f}(\mathfrak{F})
  =\prod_{i=1}^\infty p_i^{v_{\pi_2(i)}(exp(F(p_{\pi_1(i)})))}
  =\prod_{i=1}^\infty p_i^{v_{\pi_2(i)}(lcm(exp(F_1(p_{\pi_1(i)}),exp(F_2(p_{\pi_1(i)})))}\\
  =\prod_{i=1}^\infty p_i^{\max\{v_{\pi_2(i)}(exp(F_1(p_{\pi_1(i)}))), v_{\pi_2(i)}(exp(F_2(p_{\pi_1(i)})))\}}\\
  =lcm\left(\prod_{i=1}^\infty p_i^{v_{\pi_2(i)}(exp(F_1(p_{\pi_1(i)})))}, \prod_{i=1}^\infty p_i^{v_{\pi_2(i)}(exp(F_2(p_{\pi_1(i)})))}\right)=lcm(\mathfrak{f}(\mathfrak{F}_1),\mathfrak{f}(\mathfrak{F}_2)).
  \end{multline*}

$(d)$ \emph{If $\mathfrak{F}=\mathfrak{F}_1\wedge_{reg}\mathfrak{F}_2$, then $F(p)=\mathfrak{F}(gcd(exp(F_1(p)), exp(F_2(p))))$ for all $p\in\mathbb{P}$. In particular, $\mathfrak{f}(\mathfrak{F}_1\wedge_{reg}\mathfrak{F}_2)=gcd(\mathfrak{f}(\mathfrak{F}_1),\mathfrak{f}(\mathfrak{F}_2))$}.

Let $H(p)=\mathfrak{S}(gcd(exp(F_1(p)), exp(F_2(p))))$ for all $p\in\mathbb{P}$. Then formation $\mathfrak{H}$, locally defined by $H$,   is a saturated regular formation by Theorem \ref{thmc1}. Note that $F(p)\subseteq F_1(p)\cap F_2(p)$. It means that $exp(F(p))$ divides $gcd(exp(F_1(p)), exp(F_2(p)))$. According to the proof of Theorem \ref{thmc1} $\mathfrak{F}$ can be locally define by $G$ where $G(p)=\mathfrak{S}(exp(F(p)))$ for all $p\in\mathbb{P}$. Hence $G(p)\subseteq H(p)$ for all $p\in\mathbb{P}$. It means that $\mathfrak{F}\subseteq \mathfrak{H}$. Since $\mathfrak{F}$ is the greatest by inclusion element from $\mathrm{Reg}$   contained in $\mathfrak{F}_1$ and $\mathfrak{F}_2$ and $\mathfrak{H}\in\mathrm{Reg}$ is also contained in $\mathfrak{F}_1$ and $\mathfrak{F}_2$, $\mathfrak{F}=\mathfrak{H}$. Now

\begin{multline*}
  \mathfrak{f}(\mathfrak{F})
  =\prod_{i=1}^\infty p_i^{v_{\pi_2(i)}(exp(F(p_{\pi_1(i)})))}
  =\prod_{i=1}^\infty p_i^{v_{\pi_2(i)}(gcd(exp(F_1(p_{\pi_1(i)}),exp(F_2(p_{\pi_1(i)})))}\\
  =\prod_{i=1}^\infty p_i^{\min\{v_{\pi_2(i)}(exp(F_1(p_{\pi_1(i)}))), v_{\pi_2(i)}(exp(F_2(p_{\pi_1(i)})))\}}\\
  =gcd\left(\prod_{i=1}^\infty p_i^{v_{\pi_2(i)}(exp(F_1(p_{\pi_1(i)})))}, \prod_{i=1}^\infty p_i^{v_{\pi_2(i)}(exp(F_2(p_{\pi_1(i)})))}\right)=gcd(\mathfrak{f}(\mathfrak{F}_1),\mathfrak{f}(\mathfrak{F}_2)).
  \end{multline*}

$(f)$ \emph{The final step}.

From $(a)$-$(d)$ it follows that $\mathfrak{f}$ is the lattice isomorphism between $(\mathrm{Reg},\vee_{reg}, \wedge_{reg})$ and the Steinitz's  lattice $(\mathbb{SN}, lcm, gcd)$. Note that the Steinitz's  lattice is complete distributive lattice by \cite[Lemma 2]{Bezushchak2016}.

\section{Proof of Corollary \ref{cor5}}

Recall \cite{Skiba1990} (see also \cite[Corollary 6.4.5]{s9}) that a hereditary formation of soluble groups with the Shemetkov property in $\mathfrak{S}$ is saturated. Now according to \cite{Semenchuck1984} (see also \cite[Theorem 24.3]{s6} or \cite[Theorem 6.4.11]{s9}) $\mathfrak{F}\in\mathrm{Shem}$ if and only if there is a function $f:\mathbb{P}\rightarrow\mathcal{P}(\mathbb{P})$ with $p\in f(p)$ such that $\mathfrak{F}$ is locally defined by $h$ where $h(p)=\mathfrak{S}_{f(p)}$ for all $p\in\mathbb{P}$. Here $\mathcal{P}(\mathbb{P})$ denotes the power set of $\mathbb{P}$. Since every element of $\mathrm{Shem}$ contains $\mathfrak{N}$, we see that different such functions $h$ defines different formations.

Note that in this case $exp(\mathfrak{S}_{f(p)})$ is a complete supernatural number for all $p\in\mathbb{P}$. Hence $\mathfrak{f}(\mathfrak{F})$ is a complete number by its construction. From the other hand if $\omega$ is a complete number, then by $(b)$ of the proof of Theorem \ref{thmc4} there exists a function $f:\mathbb{P}\rightarrow\mathcal{P}(\mathbb{P})$ with $p\in f(p)$ such that  $\mathfrak{f}^{-1}(\omega)$ is locally defined by $h(p)=\mathfrak{S}_{f(p)}$ for all $p\in\mathbb{P}$. Thus  $\mathfrak{f}^{-1}(\omega)\in\mathrm{Shem}$.

It means that the restriction of the lattice isomorphism $\mathfrak{f}$ on $\mathrm{Shem}$ becomes bijection between $\mathrm{Shem}$ and the set of all complete supernatural numbers $\mathcal{C}$. Note that $(\mathcal{C}, lcm, gcd)$ is a  complete distributive and supplemented lattice by \cite[Lemma 3]{Bezushchak2016}. Therefore  $(\mathrm{Shem},\vee_{reg}, \wedge_{reg})$   is  a complete distributive and supplemented lattice.

\section{Proof of Corollary \ref{cor6}}

From $\mathrm{Shem}\subseteq\mathrm{Reg}$ it follows that $\mathfrak{F}_1\vee_{reg}\mathfrak{F}_2$ and $\mathfrak{F}_1\wedge_{reg}\mathfrak{F}_2$ are the least by inclusion hereditary saturated   formation with the Shemetkov property in $\mathfrak{S}$ which contains $\mathfrak{F}_1$ and $\mathfrak{F}_2$ and  the greatest by inclusion hereditary saturated   formation with the Shemetkov property in $\mathfrak{S}$ contained in $\mathfrak{F}_1$ and $\mathfrak{F}_2$   respectively for all $\mathfrak{F}_1,\mathfrak{F}_2\in\mathrm{Shem}$. Then the set $\mathrm{Shem}$  with natural partial order defines the lattice  $(\mathrm{Shem},\vee_{reg}, \wedge_{reg})$  which is  a complete distributive and supplemented lattice by Corollary \ref{cor5}. And hence has the least and the greatest elements.

{\small\bibliographystyle{spmpsci}\bibliography{Regular}}

\end{document}